\documentclass[twoside,10pt,leqno]{amsart}
\usepackage{amsfonts}
\usepackage{amsmath}
\usepackage{latexsym}
\usepackage{bbm}
\usepackage{amssymb}
\usepackage{amscd}
\numberwithin{equation}{section}

\begin{document}

\title{LARGE SIEVE INEQUALITY WITH CHARACTERS FOR POWERFUL MODULI}

\author{STEPHAN BAIER \and LIANGYI ZHAO}

\date{\today}

\maketitle

\begin{abstract}
In this paper we aim to generalize the results in \cite{Ba1},\cite{Ba2},\cite{Zha} and develop a general formula for large sieve with characters to powerful moduli that will be an improvement to the result in \cite{Zha}.\\ \\
\end{abstract}

keywords: large sieve inequality; power moduli. \\


\section{Introduction}
Throughout this paper, we reserve the symbols $c_i$ $(i=1,2,...)$ for absolute positive constants.  Large sieve was an idea originated by J. V. Linnik \cite{JVL1} in 1941 while studying the distribution of quadratic non-residues.  Refinements of this idea were made by many.  In this paper, we develop a large sieve inequality for powerful moduli.  More in particular, we aim to have an estimate for the following sum
\begin{equation}
\sum\limits_{q\le Q} \sum\limits_{\substack{ a=1\\(a,q)=1}}^{q^k}
\left\vert \sum\limits_{n= M+1}^{M+N} a_n e\left(\frac{a}{q^k}n\right)
\right\vert^2, \label{A1}
\end{equation}
where $k\ge 2$ is a natural number.  In the sequel, let
$$
Z:=\sum\limits_{n=M+1}^{M+N} \vert a_n\vert^2.
$$

With $k=1$ in \eqref{A1}, it is
\begin{equation} \label{kis1}
 \ll (Q^2+N)Z.
\end{equation}
This is in fact the consequence of a more general result first introduced by H. Davenport and H. Halberstam \cite{DH1} in which the Farey fractions in the outer sums of \eqref{A1} can be replaced by any set of well-spaced points.  Applying the said more general result, \eqref{A1} is bounded above by 
\begin{equation}
\ll (Q^{k+1}+QN)Z, \mbox{ and } \ll (Q^{2k}+N)Z\label{A4}
\end{equation}
(see \cite{Zha}).  Literature abound on the subject of the classical large sieve. See \cite{BD1}, \cite{HD}, \cite{DH1}, \cite{PXG}, \cite{JVL1}, \cite{HM}, \cite{HM2} and \cite{HLMRCV}.  In \cite{Zha} it was proved that the sum (\ref{A1}) can be estimated by 
\begin{equation}
\ll \left(Q^{k+1}+\left(NQ^{1-1/\kappa}+
N^{1-1/\kappa}Q^{1+k/\kappa}\right)N^{\varepsilon}\right)Z,\label{A5}
\end{equation}
where $\kappa:=2^{k-1}$ and the implied constant depends on $\varepsilon$ and $k$.  Furthermore, when appropriate, some of the constants $c_i$'s and the implied constants in $\ll$ in the remainder of this paper will depend on $\varepsilon$ or both $\varepsilon$ and $k$.  In \cite{Ba1} and \cite{Ba2}  this bound was improved for $k=2$.  Extending the elementary method in \cite{Ba1} 
to higher power moduli, we here establish the following bound for (\ref{A1}).\\ 

{\bf Theorem 1:} \begin{it} We have
\begin{equation} \label{A3}
 \sum\limits_{q\le Q} \sum_{\substack{a=1 \\ (a,q)=1}}^{q^k}
\left\vert \sum\limits_{n= M+1}^{M+N} a_n e\left(\frac{a}{q^k}n\right)
\right\vert^2 \ll (\log\log 10NQ)^{k+1}(Q^{k+1}+N+ N^{1/2+\varepsilon}Q^k)Z.
\end{equation}
\end{it}
For $k\ge 3$ Theorem 1 improves the classical bounds 
(\ref{A4}) as well as Zhao's bound (\ref{A5}) in the range
$N^{1/(2k)+\varepsilon}\ll Q\ll N^{(\kappa-2)/(2(k-1)\kappa-2k)-\varepsilon}$.
In particular, for $k=3$ we obtain an improvement in the range $N^{1/6+\varepsilon}\ll
Q\ll N^{1/5-\varepsilon}$. We note that for a large $k$ 
the exponent $(\kappa-2)/(2(k-1)\kappa-2k)$ is close to $1/(2(k-1))$. 

Extending the Fourier analytic 
methods in \cite{Ba2}, \cite{Zha}, we establish another bound for
cubic moduli which improves the bounds (\ref{A4}), (\ref{A5}) in the range 
$N^{7/25+\varepsilon} \ll Q \ll N^{1/3-\varepsilon}$. \\

{\bf Theorem 2:} \begin{it} Suppose that $1\le Q\le N^{1/2}$. Then we have
\begin{equation}
\sum\limits_{q\le Q} \sum_{\substack{a=1\\(a,q)=1}}^{q^3}
\left\vert \sum\limits_{n= M+1}^{M+N} a_n e\left(\frac{a}{q^3}n\right) \right\vert^2 \ll
\left\{ \begin{array}{llll} N^{\varepsilon}(Q^{4}+N^{9/10}Q^{6/5})Z, &
\mbox{ if }\ N^{7/24}\le Q\le N^{1/2},\\ \\ NQ^{6/7+\varepsilon}Z, & 
\mbox{ if }\ 1\le Q<N^{7/24}.\end{array} \right.
\end{equation}\end{it}

Unfortunately, our Fourier analytic method does not yield any improvement if $k\ge 4$.  

\section{Proof of Theorem 1}
Let ${\mathcal{S}}$ be the set of $k$-th powers of natural numbers. 
Let 
$Q_0\ge \sqrt{N}$. Set 
$$
{\mathcal{S}}(Q_0)={\mathcal{S}}\cap (Q_0,2Q_0].
$$
We first note, by classical large sieve, setting $Q=\sqrt{N}$ in \eqref{kis1},
\begin{equation}
\sum\limits_{q\le \sqrt{N}} 
\sum_{\substack{a=1\\(a,q)=1}}^q \left\vert 
\sum\limits_{n=M+1}^{M+N} a_n e\left(\frac{a}{q}n\right)\right\vert^2 \le 2NZ.
\label{C20}
\end{equation} 
Let
$$
{\mathcal{S}}_t(Q_0)=\{q\in \mathbbm{N} \ :\ tq\in {\mathcal{S}}(Q_0)\}.
$$
Let $t=p_1^{v_1}\cdots p_n^{v_n}$ be the prime decomposition of $t$.
Furthermore, let 
$$
u_i:=\left\lceil\frac{v_i}{k}\right\rceil,
$$    
where for $x\in \mathbbm{R}$, $\lceil x\rceil=\min\{k\in \mathbbm{Z}\ :\ k\ge x\}$ is the
ceiling of $x$. Moreover, set
$$
f_t=p_1^{u_1}\cdots p_n^{u_n}.
$$
Therefore, for all $q_0^k=q\in {\mathcal{S}}$, $q$ is divisible by $t$ if and only if $q_0$ is
divisible by $f_t$. Therefore, we have
$$
{\mathcal{S}}_t(Q_0)=\{q_1^k g_t\ : \ Q_0^{1/k}/f_t<q_1\le (2Q_0)^{1/k}/f_t\},
$$
where 
$$
g_t:=\frac{f_t^k}{t}.
$$
Moreover we note that
$$
{\mathcal{S}}_t(Q_0)\subset (Q_0/t,2Q_0/t]
$$
and that 
\begin{equation}
\vert {\mathcal{S}}_t(Q_0)\vert \le \frac{(2Q_0)^{1/k}}{f_t}.\label{Zu2}
\end{equation}
We set for $m\in \mathbbm{N}$, $l\in \mathbbm{Z}$ with $(m,l)=1$
\begin{equation}
A_t(u,m,l)=\max\limits_{Q_0/t\le y\le 2Q_0/t} \vert\{q\in {\mathcal{S}}_t(Q_0)\cap 
(y,y+u]\ : \
q\equiv l \mbox{ mod } m\}\vert.\label{AA}
\end{equation}
Let $\delta_t(m,l)$ be the number of solutions $x$ to the congruence
$$
x^kg_t\equiv l \mbox{ mod } m.
$$ 
We now use Theorem 2 in \cite{Ba1} with $Q_0\ge \sqrt{N}$:\\

{\bf Theorem 3:} \begin{it}  
Assume that for all $t\in \mathbbm{N}$, $m\in\mathbbm{N}$, $l\in \mathbbm{Z}$,
$u\in\mathbbm{R}$  with $t\le \sqrt{N}$, $m\le \sqrt{N}/t$, $(m,l)=1$, 
$mQ_0/\sqrt{N}\le u\le Q_0/t$ the conditions
\begin{equation}
A_t(u,m,l)\le C\left(1+\frac{\vert {\mathcal{S}}_t(Q_0)\vert /m}{Q_0/t}\cdot u\right)
\delta_t(m,l),\label{C5}
\end{equation}
\begin{equation}
\sum_{\substack{l=1\\(m,l)=1}}^{m} \delta_t(m,l)\le m, \label{C6}
\end{equation}
\begin{equation}
\delta_t(m,l)\le X\label{C7}
\end{equation} 
hold for some suitable positive numbers $C$ and $X$. Then,
\begin{equation} \label{8}
\sum\limits_{q\in {\mathcal{S}}(Q_0)} \sum_{\substack{a=1\\(a,q)=1}}^q \left\vert \sum\limits_{n=M+1}^{M+N} a_n e\left(\frac{a}{q}n\right) \right\vert^2 \le c_0C (\min\{Q_0X,N\}+Q_0) \left(\sqrt{N}\log\log 10N+ \max\limits_{r\le\sqrt{N}} \sum\limits_{t\vert r} \vert {\mathcal{S}}_t(Q_0)\vert\right)Z. 
\end{equation}
\end{it} 

First, we have to check the validity of the conditions (\ref{C5}), 
(\ref{C6}) and (\ref{C7}).
Conditions (\ref{C5}) and (\ref{C6}) are obviously satisfied with $C$ absolute. We further 
suppose that $(g_t,m)=1$ for otherwise $\delta_t(m,l)=0$ since $(m,l)=1$. 
Therefore, we must
estimate the number of solutions to
\begin{equation}
x^k\equiv \overline{g_t}l \mbox{ mod } m, \label{C8}
\end{equation}
where $\overline{g_t}$ is the multiplicative inverse of ${g_t}$ modulo $m$.
By the virtue of the Chinese remainder theorem, it suffices to estimate 
the number of solutions to (\ref{C8}) with $m$ as a prime power, say
$m=p^e$, for $p\in\mathbbm{P}$ and $e\in\mathbbm{N}$. Note that the function 
$$
\sigma_k\ : \ (\mathbbm{Z}/p^e\mathbbm{Z})^*\longrightarrow 
(\mathbbm{Z}/p^e\mathbbm{Z})^*\ : \ x \longrightarrow x^k
$$
is an endomorphism. Hence it is enough to estimate the size of its 
kernel ker$(\sigma_k)$. If $k=\pi_1^{a_1}\cdots \pi_h^{a_h}$
is the prime decomposition of $k$, then 
$$
\sigma_k=\prod\limits_{i=1}^h \sigma_{\pi_i}^{a_i}.
$$
Thus,
\begin{equation}
\vert\mbox{ker } \sigma_k\vert \le \prod\limits_{i=1}^h
\vert \mbox{ker } \sigma_{\pi_i} \vert^{a_i}.\label{Ker}
\end{equation}
Hence, it suffices to estimate the size of $\vert \mbox{ker } \sigma_{\pi}
\vert$ for prime numbers $\pi$.

For $p\in\mathbbm{P}$, 
$$
x^{\pi}-1\equiv 0\mbox{ mod } p
$$
has at most $\pi$ solutions. By elementary result 
(see \cite{Qua}, for example), 
a solution, $a$ mod $p^e$ with $e\ge 1$, of the congruence
\begin{equation}
x^{\pi}-1 \equiv 0 \mbox{ mod } p^e \label{C10}
\end{equation}
lifts to more than one solution to
$$
x^{\pi}-1 \equiv 0 \mbox{ mod } p^{e+1}
$$  
only when $p\vert \pi a^{\pi-1}$ and $p^{e+1}\vert a^{\pi}-1$. If $p\not=\pi$, 
$p\vert \pi a^{\pi-1}$ implies $p\vert a$, but it is not possible that
$p^{e+1}\vert a^{\pi}-1$ as $(a^{\pi}-1,a)=1$. Thus, in this case 
(\ref{C10}) has at most $\pi$
solutions for all $e\ge 1$. In the following, we consider the case $p=\pi$. 

By Fermat's little theorem, there exists only one solution of the 
congruence
$$
x^{\pi}-1\equiv 0 \mbox{ mod } \pi,
$$
namely $1$ mod $\pi$. This solution lifts to exactly $\pi$ solutions to
$$
x^{\pi}-1\equiv 0 \mbox{ mod } \pi^2,
$$
namely
$$
1,\ 1+\pi,\ 1+2\pi,\ ...,1+(\pi-1)\pi \mbox{ mod } \pi^2.
$$
More generally, if $a$ mod $\pi^e$ is a solution to 
\begin{equation}
x^{\pi}-1 \equiv 0 \mbox{ mod } \pi^e,\label{CO1}
\end{equation}
then, if $a$ lifts to solutions to
$$
x^{\pi}-1 \equiv 0 \mbox{ mod } \pi^{e+1},
$$
they are of the form
\begin{equation}
a,\ a+\pi^e,\ a+2\pi^e,\ ...,\ a+(\pi-1)\pi^e \mbox{ mod } \pi^{e+1}.
\label{CO2}
\end{equation}
Assume there are $j_1,j_2\in\{0,...,\pi-1\}$, $j_1\not= j_2$ such that both 
$a+j_1\pi^e$ and $a+j_2\pi^e$ lift to solutions modulo $\pi^{e+2}$. Then
$\pi^{e+2} \vert (a+j_1\pi^e)^{\pi}-1$ and 
$\pi^{e+2} \vert (a+j_2\pi^e)^{\pi}-1$, hence
$$
(a+j_1\pi^e)^{\pi}-(a+j_2\pi^e)^{\pi}
= (j_1-j_2)\pi^e\sum\limits_{i=0}^{\pi-1}
(a+j_1\pi^e)^{\pi-1-i}(a+j_2\pi^e)^{i}
$$
is divisible by $\pi^{e+2}$. If $e\ge 2$, this implies $a\equiv 0$ mod $\pi$, 
but then $a$ cannot be a solution to (\ref{CO1}). Therefore, if $e\ge 2$, only
one of the solutions (\ref{CO2}) lifts to a solution modulo $\pi^{e+2}$.
From this we infer that the 
number of solutions to (\ref{CO1}) never exceeds $\pi^2$, {\it i.e.}
$$
\vert \mbox{ker } \sigma_{\pi}\vert \le \pi^2.
$$
Combining this with (\ref{Ker}), we get
$$
\vert \mbox{ker } \sigma_{k}\vert \le k^2.
$$
Therefore, by the Chinese remainder theorem, we obtain    
$$
\delta_t(m,l)\le k^{2\omega(m)},
$$
where $\omega(m)$ is the number of distinct prime divisors of $m$.
Since $2^{\omega(m)}$ is the number of square-free divisors of $m$,
we have 
$$
k^{2\omega(m)}\le \tau(m)^{2\log_2 k}\ll m^{\varepsilon},
$$
where $\tau(m)$ is the number of divisors of $m$.
Thus, if $m\le \sqrt{N}$, (\ref{C7}) holds with
$$
X\ll N^{\varepsilon}.
$$

Now, by Theorem 3,
\begin{equation}
\sum\limits_{q\in {\mathcal{S}}(Q_0)} 
\sum\limits_{\substack{ a=1\\(a,q)=1}}^q \left\vert 
\sum\limits_{n=M+1}^{M+N} a_n e\left(\frac{a}{q}n\right)\right\vert^2
\label{C12}
\end{equation} 
is majorized by
$$
\ll (\min\{Q_0N^{\varepsilon},N\}+Q_0)\left(\sqrt{N}\log\log(10N)+
\max\limits_{r\le \sqrt{N}} \sum\limits_{t\vert r} Q_0^{1/k}f_t^{-1}\right)Z.
$$
The function 
$$
G(r)=\sum\limits_{t\vert r} \frac{1}{f_t}
$$
is clearly multiplicative. If $r$ is a prime power $p^v$, then
$$
G(r)\le 1+k\left(\frac{1}{p}+\frac{1}{p^2}+...\right)=1+\frac{k}{p-1}
\le \left(1+\frac{1}{p-1}\right)^k=\left(\frac{p^v}{\varphi\left(p^v\right)}
\right)^k.
$$
Hence, for all $r\in \mathbbm{N}$ we have
\begin{equation}
G(r)\le \left(\frac{r}{\varphi(r)}\right)^k\ll (\log \log 10r)^k.
\label{Zu1}
\end{equation}
Hence (\ref{C12}) is
$$
\ll (\log\log 10NQ_0)^{k+1}(\sqrt{N}+Q_0^{1/k})(
\min\{Q_0N^{\varepsilon},N\}+Q_0).
$$
The above is always majorized by
$$
\ll (\log\log 10NQ_0)^{k+1}\left(Q_0^{1+1/k}+N^{1/2+\varepsilon}Q_0\right).
$$
Summing over all relevant dyadic intervals and combining with (\ref{C20}), we see that (\ref{A1}) is majorized by
$$
\ll (\log\log 10NQ)^{k+1}(Q^{k+1}+N+N^{1/2+\varepsilon}Q^k)Z.
$$
Therefore, our result follows. $\Box$ 

\section{Proof of Theorem 2}
\subsection{Reduction to Farey fractions in short intervals}
As in \cite{Ba1}, \cite{Ba2}, our starting point is the following general large sieve inequality.\\  
  
{\bf Lemma 1:} \begin{it} Let $\left(\alpha_r\right)_{r\in\mathbbm{N}}$ be a
sequence of real numbers. Suppose that $0<\Delta\le 1/2$ and 
$R\in \mathbbm{N}$. Put 
\begin{equation}
K(\Delta):=\max\limits_{\alpha\in \mathbbm{R}} 
\sum_{\substack{ r=1\\ \| \alpha_r -\alpha\|\le \Delta }}^R 1,\label{1}
\end{equation}
where $\| x \|$ denotes the distance of a real $x$
to its closest integer.
Then 
$$
\sum\limits_{r=1}^R \left\vert S\left(\alpha_r\right)\right\vert^2
\le c_{1}K(\Delta)(N+\Delta^{-1})Z.
$$
\end{it}

In the sequel, we suppose that ${\mathcal{S}}$ is the set of 
cubes of natural numbers and that    
$\alpha_1,...,\alpha_R$ is the sequence
of Farey fractions $a/q$ with $q\in {\mathcal{S}}(Q_0)$,  
$1\le a\le q$ and $(a,q)=1$, where $Q_0\ge 1$.  
We further suppose that $\alpha\in \mathbbm{R}$ and $0<\Delta\le 1/2$. 
Put
$$
I(\alpha):=[\alpha-\Delta,\alpha+\Delta]
\ \ \mbox{ and }\ \ 
P(\alpha):= \sum_{\substack{ q\in {\mathcal{S}}\cap(Q_0,2Q_0]\\ (a,q)=1 \\ a/q\in I(\alpha)}} 1.
$$
Then we have
$$
K(\Delta)=\max\limits_{\alpha\in \mathbbm{R}} P(\alpha).\label{60}
$$
Therefore, the proof of Theorem 
2 reduces to estimating $P(\alpha)$.   

As in \cite{Ba1} and \cite{Ba2}, we begin with an idea of D. Wolke \cite{Wol}. Let $\tau$ be a positive number with
\begin{equation} \label{P1}
1\le \tau\le \frac{1}{\sqrt{\Delta}}.
\end{equation}
In \cite{Ba1} and \cite{Ba2} we put
$\tau:=1/\sqrt{\Delta}$, but in fact our method works for all $\tau$
satisfying $(\ref{P1})$. We will later fix $\tau$ in an optimal manner. 
In the said earlier papers, $\tau=1/\sqrt{\Delta}$ was the optimal choice. 

By Dirichlet's approximation theorem, $\alpha$ can be written in the form
$$
\alpha=\frac{b}{r}+z,
$$
where 
\begin{equation}
r\le \tau,\ (b,r)=1,\ 
\vert z\vert \le \frac{1}{r\tau}.\label{P2}
\end{equation}
Thus, it suffices to estimate $P(b/r+z)$ for all $b,r,z$ satisfying
(\ref{P2}).

We further note that we can restrict ourselves to the case when 
\begin{equation}
z\ge \Delta.\label{P4}
\end{equation}
If $\vert z\vert<\Delta$, then
$$ 
P(\alpha)\le P\left(\frac{b}{r}-\Delta\right)+P\left(\frac{b}{r}+\Delta\right).
$$ 
Furthermore, we have
$$
\Delta\le \frac{1}{\tau^2}\le \frac{1}{r\tau}.
$$ 
Therefore this case can be reduced to the case $\vert z\vert=\Delta$. 
Moreover, as $P(\alpha)=P(-\alpha)$, we can choose
$z$ positive. So we can assume (\ref{P4}), without any loss of generality.

Summarizing the above observations, we deduce\\

{\bf Lemma 2:} \begin{it} We have
\begin{equation}
K(\Delta)\le 2\max_{\substack{ r\in \mathbbm{N} \\ r\le \tau}} \max_{\substack{b\in \mathbbm{Z}\\ (b,r)=1}} \max_{\Delta\le z\le 1/(\tau r)} P\left(\frac{b}{r}+z\right).
\end{equation} \end{it}

\subsection{Estimation of $P(b/r+z)$ - first way} 
We now prove a first estimate for $P\left(b/r+z\right)$ by using some results 
in \cite{Ba1}.  In the sequel, we suppose that the conditions (\ref{P1}), 
(\ref{P2}) and (\ref{P4}) are satisfied.

By inequality (41) in \cite{Ba1}, we have 
\begin{equation}
P\left(\frac{b}{r}+z\right)\le 1+6 \sum\limits_{t\vert r} \sum_{\substack{0<m\le 4rzQ_0/t\\ (m,r/t)=1}}
A_t\left(\frac{\Delta Q_0}{tz},\frac{r}{t},-\overline{b}m\right),\label{P14}
\end{equation}
where $A_t(u,m,l)$ is defined as in (\ref{AA}) and $\overline{b}$ is the multiplicative inverse of $b$ modulo $r$.  By the results of section 2, for ${\mathcal{S}}$ the set of cubes, 
the conditions (\ref{C5}), (\ref{C6}) and (\ref{C7}) with 
$X=\Delta^{-\varepsilon}$ are
satisfied for all
$t\in \mathbbm{N}$, $m\in\mathbbm{N}$, $l\in \mathbbm{Z}$,
$u\in\mathbbm{R}$  with $t\le \tau$, $m\le \tau/t$, 
$(m,l)=1$, 
$mQ_0/\tau\le u\le Q_0/t$.
Conditions (\ref{C5}) and (\ref{C7}) imply
\begin{equation}\label{AA1}
\sum_{\substack{0<m\le 4rzQ_0/t \\ (m,r/t)=1}} A_t\left(\frac{\Delta Q_0}{tz},\frac{r}{t},-\overline{b} m\right) \le C\left(1+\frac{\Delta t\vert {\mathcal{S}}_t(Q_0)\vert}{rz}\right) \frac{4rzQ_0X}{t} 
\end{equation}
From (\ref{P14}), (\ref{AA1}) and 
$$
\sum\limits_{t\vert r}\frac{1}{t}\le \prod\limits_{p\vert r}
\left(1+\frac{1}{p}+\frac{1}{p^2}+...\right) = \prod\limits_{p\vert r} 
\frac{p}{p-1} =
\frac{r}{\varphi(r)}\le c_{2}\log\log 10r,
$$
we derive
\begin{equation}
P\left(\frac{b}{r}+z\right)\le 1+
c_{3}Q_0X\left(rz\log\log 10r +\Delta \sum\limits_{t\vert r} \vert 
{\mathcal{S}}_t(Q_0)\vert\right).
\label{AA5}
\end{equation} 
Furthermore, by (\ref{Zu2}) and (\ref{Zu1}), we have 
$$
\sum\limits_{t\vert r} \vert {\mathcal{S}}(Q_0) \vert \ll 
(\log\log 10r)^3 Q_0^{1/3}.
$$
Thus, from (\ref{AA5}) and the fact that $r \leq \tau = \Delta^{-1/2}$, we obtain\\

{\bf Proposition 1:} \begin{it} Let ${\mathcal{S}}$ be the set of cubes of
natural numbers. 
Suppose that
the conditions (\ref{P1}), 
(\ref{P2}) and (\ref{P4}) are satisfied. Then we have    
\begin{equation}
P\left(\frac{b}{r}+z\right)\le 1+
c_{4}\Delta^{-\varepsilon}\left(Q_0^{4/3}\Delta+Q_0rz\right).\label{AA6}
\end{equation} \end{it}

\subsection{Estimation of $P(b/r+z)$ - second way}
We now prove a second estimate for $P\left(b/r+z\right)$ by extending the 
Fourier ana\-lytic methods in
\cite{Ba2}, \cite{Zha} to cubic moduli. 
The following bound for $P(b/r+z)$ can be proved in the same manner as Lemma 2 
in [2].\\
       
{\bf Lemma 3:} \begin{it}
Let ${\mathcal{S}}$ be the set of cubes of natural numbers. Suppose that 
\begin{equation}
\frac{Q_0\Delta}{z}\le \delta\le Q_0. \label{P10}
\end{equation}
Then,
\begin{equation} \label{PP}
P\left(\frac{b}{r}+z\right) \le c_5\left(1+\frac{1}{\delta} \int\limits_{Q_0}^{2Q_0} \Pi(\delta, y) \,{\rm d}y\right),
\end{equation}
where $I(\delta, y)= [y^{1/3}-c_6\delta/Q_0^{2/3},  y^{1/3}+c_6\delta/Q_0^{2/3}]$, $J(\delta, y) = [(y-4\delta)rz, (y+4\delta)rz]$ and
\begin{equation} \label{defPi}
\Pi(\delta, y) = \sum\limits_{q \in I(\delta, y)} \sum_{\substack{m \in J(\delta, y) \\ m \equiv -bq^3 \mod{r} \\ m\not=0 }} 1 .
\end{equation}
\end{it}

We shall prove the following\\

{\bf Proposition 2:} \begin{it} 
Let ${\mathcal{S}}$ be the set of cubes of natural numbers. Suppose that
the conditions (\ref{P1}), 
(\ref{P2}) and (\ref{P4}) are satisfied. Then we have   
\begin{equation}
P\left(\frac{b}{r}+z\right)\le 
c_7\Delta^{-\varepsilon}\left(Q_0^{4/3}\Delta
+Q_0^{1/3}\Delta r^{-1/3}z^{-1}+\Delta^{-1/2}(rz)^{1/2}\right).\label{QQ}
\end{equation}\end{it}

To derive Proposition 2 from Lemma 3, we need the following standard results 
from 
Fourier analysis.\\ 

{\bf Lemma 4:} (Poisson summation formula, \cite{Bum}) \begin{it}  
 Let $f(x)$ be a complex-valued 
function on the real numbers that is piecewise continuous with only finitely 
many discontinuities and for all real numbers $a$ satisfies
$$
f(a)=\frac{1}{2}\left(\lim\limits_{x\rightarrow a^-} f(x) +
\lim\limits_{x\rightarrow a^+} f(x)\right).
$$
Moreover, suppose that $f(x)\le c_8(1+\vert x\vert)^{-c}$ for some $c>1$.
Then,
$$
\sum\limits_{n\in \mathbbm{Z}} f(n) = \sum\limits_{n\in \mathbbm{Z}} 
\hat{f}(n), \; \mbox{where} \; \hat{f}(x):=\int\limits_{-\infty}^{\infty} f(y)e(xy) {\rm d}y,
$$
the Fourier transform of $f(x)$. \end{it}\\

{\bf Lemma 5:} (see \cite{Zha}, for example) \begin{it}
For $x\in \mathbbm{R}\setminus \{0\}$ define
$$
\phi(x):=\left(\frac{\sin \pi x}{2x}\right)^2, \; \mbox{and} \; \phi(0):=\lim\limits_{x\rightarrow 0}\phi(x)=\frac{\pi^2}{4}. 
$$
Then $\phi(x)\ge 1$ for $\vert x\vert \le 1/2$, and
the Fourier transform of the function $\phi(x)$ is
$$
\hat{\phi}(s)=\frac{\pi^2}{4}\max\{1-\vert s\vert, 0\}.
$$\end{it}

{\bf Lemma 6:} (see Lemma 3.1. in \cite{Kol}) \begin{it} Let $F$ $:$ 
$[a,b]\rightarrow \mathbbm{R}$ be twice differentiable. 
Assume that $\vert F^{\prime}(x)\vert \ge u>0$ 
for all $x\in [a,b]$. Then,
$$
\left\vert \int\limits_{a}^{b} e^{iF(x)} {\rm d}x\right\vert \le
\frac{c_9}{u}.
$$
\end{it}

{\bf Lemma 7:} (see Lemma 4.3.1. in \cite{Bru}) \begin{it} Let $F$ $:$ 
$[a,b]\rightarrow \mathbbm{R}$ be twice continuously differentiable. 
Assume that $\vert F^{\prime\prime}(x)\vert \ge u>0$ 
for all $x\in [a,b]$. Then,
$$
\left\vert \int\limits_{a}^{b} e^{iF(x)} {\rm d}x\right\vert \le
\frac{c_{10}}{\sqrt{u}}.
$$\end{it}

We shall also need the following estimates for cubic exponential sums.\\

{\bf Lemma 8:} (see \cite{Hua}, \cite{Vin}) \begin{it}
Let $c\in \mathbbm{N}$, $k,l\in 
\mathbbm{Z}$ with $(k,c)=1$. Then,
$$
\sum\limits_{d=1}^c e\left(\frac{kd^3+ld}{c}\right) \le 
c_{11}c^{1/2+\varepsilon}(l,c).
$$ 
Furthermore,
$$
\sum\limits_{d=1}^c e\left(\frac{kd^3}{c}\right) \le 
c_{11}c^{2/3}.
$$
\end{it}\\

{\bf Proof of Proposition 2:} We put 
\begin{equation}
\delta:=\frac{Q_0\Delta}{z}.\label{Delta}
\end{equation}
By Lemma 5, \eqref{defPi} can be estimated by
\begin{equation} \label{P13}
\Pi (\delta, y) \le \sum\limits_{q\in \mathbbm{Z}} \ \phi\left(\frac{q-y^{1/3}}{2c_6\delta/Q_0^{2/3}} \right) \sum_{\substack{m\in \mathbbm{Z} \\ m\equiv -bq^3 \mod{r}}} \phi\left(\frac{m-yrz}{8\delta rz}\right).
\end{equation}

Using Lemma 4 after a linear change of variables, we transform the inner sum 
on the right-hand side of (\ref{P13}) into 
$$
\sum_{\substack{m\in \mathbbm{Z} \\ m\equiv -bq^3 \mod{r} }} \phi\left(\frac{m-yrz}{8\delta rz}\right) = 8\delta z \sum\limits_{j\in \mathbbm{Z}} e\left(\frac{jbq^3}{r}+jyz\right) \hat{\phi}(8j\delta z).
$$ 
Therefore, we get for the double sum on the right-hand side of (\ref{P13})
\begin{equation} \label{F1}
\begin{split}
\sum\limits_{q\in \mathbbm{Z}} \ &
\phi\left(\frac{q-y^{1/3}}{2c_6\delta/Q_0^{2/3}} 
\right) \sum_{\substack{ m\in \mathbbm{Z}\smallskip\\ m\equiv -bq^3 \mod{r}}} \phi\left(\frac{m-yrz}{8\delta rz}\right) \\
&= 8\delta z \sum_{j\in \mathbbm{Z}} e(jyz)\hat{\phi}(8j\delta z) 
\sum_{d=1}^{\tilde r} e\left(\frac{\tilde jbd^3}{\tilde r}\right)
\sum_{\substack{ k\in \mathbbm{Z}\\ k\equiv d \mod{\tilde r}}} \phi\left(\frac{k-y^{1/3}}{2c_6\delta/Q_0^{2/3}}\right),
\end{split}
\end{equation}
where $\tilde r:=r/(r,j)$ and $\tilde j:=j/(r,j)$. 
Again using Lemma 4 after a linear 
change of variables, we transform the inner sum 
on the right-hand side of (\ref{F1}) into
\begin{equation}
\sum\limits_{\substack{ k\in \mathbbm{Z}\\ k\equiv d \mod{\tilde r}}}
\phi\left(\frac{k-y^{1/3}}{2c_6\delta/Q_0^{2/3}}\right)=
\frac{2c_6\delta}{\tilde r Q_0^{2/3}} \sum\limits_{l\in \mathbbm{Z}} 
e\left(l\cdot\frac{d-y^{1/3}}{\tilde r}\right)\hat{\phi}
\left(\frac{2c_6l\delta}{\tilde r Q_0^{2/3}}\right).\label{F2}
\end{equation}
From (\ref{F1}) and (\ref{F2}), we obtain 
\begin{equation} \label{F3}
\begin{split}
\frac{1}{\delta} & \int\limits_{Q_0}^{2Q_0} \
\sum\limits_{q\in \mathbbm{Z}} \ 
\phi\left(\frac{q-y^{1/3}}{2c_6\delta/Q_0^{2/3}} 
\right) \sum_{\substack{ m\in \mathbbm{Z}\\ m\equiv -bq^3 \mod{r}}} \phi\left(\frac{m-yrz}{8\delta rz}\right) {\rm d}y \\
&\le \frac{16c_6\delta z}{Q_0^{2/3}} \sum\limits_{j\in \mathbbm{Z}} 
\frac{\hat{\phi}(8j\delta z)}{\tilde r} \sum\limits_{l\in \mathbbm{Z}}
\hat{\phi}\left(\frac{2c_6l\delta}{\tilde r Q_0^{2/3}}\right)\left\vert 
\sum\limits_{d=1}^{\tilde r} 
e\left(\frac{\tilde j bd^3+ld}{\tilde r}\right) \int\limits_{Q_0}^{2Q_0} e\left(jyz-l\cdot\frac{y^{1/3}}{\tilde r}\right)\
{\rm d}y\right\vert.
\end{split}
\end{equation}
Applying the Lemmas 5 and 8 to the right-hand side of (\ref{F3}), and 
taking $r\le 1/\sqrt{\Delta}$ by (\ref{P1}) and (\ref{P2}) into account, we 
deduce
\begin{equation} \label{F4}
\begin{split}
\frac{1}{\delta} & \int\limits_{Q_0}^{2Q_0} \ \sum\limits_{q\in \mathbbm{Z}} \ \phi\left(\frac{q-y^{1/3}}{c_6\delta/Q_0^{2/3}} \right) \sum_{\substack{ m\in \mathbbm{Z}\\ m\equiv -bq^3 \mod{r}}} \phi\left(\frac{m-yrz}{8\delta rz}\right) {\rm d}y \\
&\le \frac{c_{12}\delta z\Delta^{-\varepsilon}}{Q_0^{2/3}} \left( \sum\limits_{\vert j \vert \le 1/(8\delta z)} \frac{1}{\sqrt{\tilde r}} \sum_{\substack{ \vert l\vert \le (\tilde r Q_0^{2/3})/(2c_6\delta)\\ l\not= 0 }} (l,\tilde r) \left\vert \int\limits_{Q_0}^{2Q_0} e\left(jyz-l\cdot\frac{y^{1/3}}{\tilde r}\right)\ {\rm d}y \right\vert+\sum\limits_{\vert j \vert \le 1/(8\delta z)} \frac{1}{\sqrt[3]{\tilde r}} \left\vert \int\limits_{Q_0}^{2Q_0} e\left(jyz\right)\ {\rm d}y \right\vert\right).
\end{split}
\end{equation}

If $j\not= 0$, then 
$$
\left\vert 
\int\limits_{Q_0}^{2Q_0} e\left(jyz\right)\ 
{\rm d}y
\right\vert \le \frac{1}{\vert j\vert z}.
$$
If $j=0$ and $l\not=0$, then 
$$
\left\vert 
\int\limits_{Q_0}^{2Q_0} e\left(jyz-l\cdot\frac{y^{1/3}}{\tilde r}\right)
\ {\rm d}y\right\vert \le \frac{c_{13}Q_0^{2/3}}{\vert l\vert}
$$
by Lemma 6 (take into account that $\tilde r=1$ if $j=0$).
If $j\not=0$ and $l\not=0$, then Lemma 7 yields
$$
\left\vert 
\int\limits_{Q_0}^{2Q_0} e\left(jyz-l\cdot\frac{y^{1/3}}{\tilde r}\right)
\ {\rm d}y\right\vert \le 
\frac{c_{14}\sqrt{\tilde r}Q_0^{5/6}}{\sqrt{\vert l\vert}}.
$$
Therefore, the right-hand side of (\ref{F4}) is majorized by
\begin{equation} \label{F5}
\le c_{15}\delta \Delta^{-\varepsilon}\left(zQ_0^{1/3}+\frac{1}{Q_0^{2/3}} \sum\limits_{1\le j \le 1/(8\delta z)} \frac{1}{j\sqrt[3]{\tilde r}}\right. +z\sum\limits_{1\le l\le Q_0^{2/3}/(2c_6\delta)} \frac{1}{l} +\left. zQ_0^{1/6}\sum\limits_{1\le j \le 1/(8\delta z)}\ \sum\limits_{1\le l\le \tilde r Q_0^{2/3}/(2c_6\delta)}\frac{(l,\tilde r)}{\sqrt{l}}\right).
\end{equation} 

Now, we estimate the sums in the last line of (\ref{F5}). Using  
(\ref{P1}), (\ref{P2}) and (\ref{Delta}), we obtain
\begin{equation}
\sum\limits_{1\le l\le Q_0^{2/3}/(2c_6\delta)} \frac{1}{l} \le 
c_{16}\Delta^{-\varepsilon}.\label{FF0}
\end{equation}
Using the definition of $\tilde r$, ({\ref{P1}), (\ref{P2}) and  (\ref{Delta}), we
obtain
\begin{equation} \label{F6}
\sum\limits_{1\le j \le 1/(8\delta z)} \frac{1}{j\sqrt[3]{\tilde r}}= \frac{1}{\sqrt[3]{r}} \sum\limits_{t\vert r} \sqrt[3]{t} \sum\limits_{\substack{ 1\le j \le 1/(8\delta z) \\ (r,j)=t }} \frac{1}{j} \le \frac{c_{17}\Delta^{-\varepsilon}}{\sqrt[3]{r}} \sum\limits_{t\vert r} t^{-2/3} \le c_{18}\Delta^{-2\varepsilon}r^{-1/3}.
\end{equation}
For $A\ge 1$, we have
\[  \sum\limits_{1\le l\le A} \frac{(l,\tilde r)}{\sqrt{l}} \le \sum\limits_{t\vert \tilde r} t \sum\limits_{1\le l\le A/t} \frac{1}{\sqrt{lt}} \ll \sqrt{A}\sum\limits_{t\vert \tilde r} 1 \ll {\tilde r}^{\varepsilon} \sqrt{A}. \]
Therefore,
\begin{equation}
\sum\limits_{1\le j \le 1/(8\delta z)}\sum\limits_{1\le l\le \tilde 
r Q_0^{2/3}/(2c_6\delta)}\frac{(l,\tilde r)}{\sqrt{l}} \le
\frac{c_{19} \Delta^{-\varepsilon}Q_0^{1/3}}{\sqrt{\delta}}  
\sum\limits_{1\le j \le 1/(8\delta z)} \sqrt{\tilde r}.\label{FF}
\end{equation}
Using the definition of $\tilde r$, we obtain 
\begin{equation} \label{F7} 
\sum\limits_{1\le j \le 1/(8\delta z)} \sqrt{\tilde r} = {\sqrt{r}} \sum\limits_{t\vert r} \frac{1}{\sqrt{t}} \sum\limits_{\substack{ 1\le j \le 1/(8\delta z) \\ (r,j)=t }} 1 \le \frac{\sqrt{r}}{8\delta z} \sum\limits_{t\vert r} \frac{1}{t^{3/2}} \le \frac{c_{20}\sqrt{r}}{\delta z}.
\end{equation}

Combining Lemma 3 and (\ref{F4}-\ref{F7}), we obtain
\begin{equation}
P\left(\frac{b}{r}+z\right)\le 
c_7\Delta^{-3\varepsilon}\left(1+\delta zQ_0^{1/3}
+\delta Q_0^{-2/3}r^{-1/3}+
\delta^{-1/2}Q_0^{1/2}\sqrt{r}\right).\label{Q1}
\end{equation} 
From (\ref{Delta}) and (\ref{Q1}), we infer the desired estimate.  Note that the first term in the right-hand side of \eqref{Q1} can be absorbed into the last term on the right-hand side of \eqref{QQ} by \eqref{P4}. $\Box$
 
\subsection{Final proof of Theorem 2}
Combining Propositions 1,2 and (\ref{P2}), we obtain  
\begin{equation}
P\left(\frac{b}{r}+z\right)\le 
c_{21}\Delta^{-\varepsilon}\left(Q_0^{4/3}\Delta+
\min\left\{Q_0rz,Q_0^{1/3}\Delta r^{-1/3}z^{-1}\right\}+
\Delta^{-1/2}\tau^{-1/2}\right).\label{E1}
\end{equation} 
If
$$
z\le \Delta^{1/2}Q_0^{-1/3}r^{-2/3},
$$
then 
$$
\min\left\{Q_0rz,Q_0^{1/3}\Delta r^{-1/3}z^{-1}\right\}=
Q_0rz\le Q_0^{2/3}\Delta^{1/2}r^{1/3}.
$$
If 
$$
z> \Delta^{1/2}Q_0^{-1/3}r^{-2/3},
$$
then 
$$
\min\left\{Q_0rz,Q_0^{1/3}\Delta r^{-1/3}z^{-1}\right\}=
Q_0^{1/3}\Delta r^{-1/3}z^{-1}\le Q_0^{2/3}\Delta^{1/2}r^{1/3}.
$$
From the above inequalities and (\ref{P2}), we deduce
\begin{equation}
\min\left\{Q_0rz,Q_0^{1/3}\Delta r^{-1/2}z^{-1}\right\}\le 
Q_0^{2/3}\Delta^{1/2}r^{1/3}\le Q_0^{2/3}\Delta^{1/2}\tau^{1/3}. \label{E2}
\end{equation}
Combining (\ref{E1}) and (\ref{E2}), we get
\begin{equation}
P\left(\frac{b}{r}+z\right)\le 
c_{22}\Delta^{-\varepsilon}\left(Q_0^{4/3}\Delta \tau^{\varepsilon}+Q_0^{2/3}\Delta^{1/2}
\tau^{1/3+\varepsilon}+
\Delta^{-1/2}\tau^{-1/2}\right).
\label{E4}
\end{equation}
Now we choose 
$$
\tau:=\left\{ \begin{array}{llll} N^{6/5}Q_0^{-4/5}, &
\mbox{ if }\ N^{7/8}\le Q_0\le N^{3/2},\\ \\ Q_0^{4/7}, & \mbox{ if }\
1 \le Q_0<N^{7/8},\end{array} \right. \; \mbox{and} \; \Delta:=\left\{ \begin{array}{llll} N^{-1}, &
\mbox{ if }\ N^{7/8}\le Q_0\le N^{3/2},\\ \\ Q_0^{-8/7}, & \mbox{ if }\
1\le Q_0<N^{7/8}.\end{array} \right. 
$$
Then the condition \eqref{P1} is satisfied in each case, and from (\ref{E4}) and Lemmas 1,2, we obtain 
\begin{equation} \label{E5}
\sum\limits_{Q_0^{1/3}\le q\le (2Q_0)^{1/3}} \sum\limits_{\substack{ a=1\\(a,q)=1}}^{q^3} \left\vert S\left(\frac{a}{q^3}\right) \right\vert^2 \ll \left\{ \begin{array}{llll} N^{\varepsilon}\left(Q_0^{4/3}+N^{9/10}Q_0^{2/5}\right)Z, & 
\mbox{ if }\ N^{7/8}\le Q_0\le N^{3/2},\\ \\ NQ_0^{2/7+\varepsilon}Z, & \mbox{ if }\
1\le Q_0<N^{7/8}.\end{array} \right.
\end{equation}
We can divide the interval $[1,Q]$ into $O(\log Q)$ subintervals of the
form $\left[Q_0^{1/3},(2Q_0)^{1/3}\right]$, where $1\le Q_0\le Q^3$. Hence,
the result of Theorem 2 follows from (\ref{E5}). $\Box$
 
\section*{Acknowledgments}
This paper was written when the first-named author held a postdoctoral position at the Harish-Chandra Research Institute at Allahabad (India) and the second-named author was supported by a postdoctoral fellowship at the University of Toronto.  The authors wish to thank these institutions for their financial support.

\vspace*{1cm}
\noindent Stephan Baier \newline
Harish-Chandra Research Institute, Chhatnag Road, Jhusi, Allahabad 211 019, India \newline
Email: {\tt sbaier@mri.ernet.in} \newline

\noindent Liangyi Zhao \newline
Department of Mathematics, University of Toronto, 100 Saint George Street, Toronto, ON M5S 3G3, Canada \newline
Email: {\tt lzhao@math.toronto.edu}


\begin{thebibliography}{9} 
\bibitem[1]{Ba1} S. Baier, {\it On the large sieve with a sparse set of moduli}, preprint.
\bibitem[2]{Ba2} S. Baier, {\it The large sieve with square moduli}, preprint. 
\bibitem[3]{BD1} E. Bombieri and H. Davenport, {\it Some inequalities involving trigonometrical polynomials}, Annali Scuola Normale Superiore - Pisa 23 (1969) 223-241.
\bibitem[4]{Bru} J. Br\"udern, {\it Einf\"uhrung in die analytische Zahlentheorie}, Springer-Verlag, Berlin ect., 1995.
\bibitem[5]{Bum} D. Bump, {\it Automorphic Forms and Representations}, Cambridge Stud. Adv. Math. 55, Cambridge Univ. Press, Cambridge, 1996.
\bibitem[6]{HD} H. Davenport, {\it Multiplicative Number Theory}, Third Edition, Graduate Texts in Mathematics, 74, Springer-Verlag, New York, etc., 2000.
\bibitem[7]{DH1} H. Davenport and H. Halberstam, {\it The values of a trigonometric polynomial at well spaced points}, Mathematika 13 (1966) 91-96, {\it Corrigendum and addendum}, Mathematika 14 (1967) 299-232.
\bibitem[8]{PXG} P. X. Gallagher, {\it The large sieve}, Mathematika 14 (1967) 14-20.
\bibitem[9]{Kol} S.W. Graham, G. Kolesnik, {\it Van der Corput's Method of Exponential Sums}, Cambridge University Press, Cambridge ect., 1991.
\bibitem[10]{JVL1} J. V. Linnik, {\it The large sieve}, Doklady Akad. nauk SSSR 36 (1941) 119-120.
\bibitem[11]{Hua} L. K. Hua, {\it On exponential sums}, Sci. Record 1 (1957) 1-4.
\bibitem[12]{HM} H. L. Montgomery, {\it Topics in Multiplicative Number Theory}, Lecture Notes in Mathematics, 227, Spring-Verlag, Berlin, etc., 1971.
\bibitem[13]{HM2} H. L. Montgomery, {\it The analytic principles of large sieve}, Bull. Amer. Math. Soc., 84 (1978) 547-567.
\bibitem[14]{HLMRCV} H. L. Montgomery and R. C. Vaughan, {\it The Large Sieve}, Mathematika 20 (1973) 119-134.
\bibitem[15]{Qua} I. Niven, H.S. Zuckerman, H.L. Montgomery, {\it An introduction to the Theory of Numbers}, John Wiley \& Sons, 
New York, 1991.  
\bibitem[16]{Vau} R.C. Vaughan, {\it The Hardy-Littlewood Method},  
Cambridge University Press, Cambridge, 1997.
\bibitem[17]{Vin} I.M. Vinogradov, {\it The Method of Trigonometrical Sums in the Theory of Numbers}, Trav. Inst. Math. Stekloff  23, 1947.
\bibitem[18]{Wol} D. Wolke, {\it On the large sieve with primes}, Acta Math. Acad. Sci. Hungar.  22  (1971/72) 239-247. 
\bibitem[19]{Zha}  L. Zhao, {\it Large sieve inequality with characters to square moduli}, Acta Arith. 112 (2004) 297-308.
\end{thebibliography}
\end{document}